\documentclass[11pt]{amsart}

\theoremstyle{plain}
\newtheorem{theorem}                {Theorem}      [section]
\newtheorem{proposition}  [theorem]  {Proposition}
\newtheorem{corollary}    [theorem]  {Corollary}
\newtheorem{lemma}        [theorem]  {Lemma}

\theoremstyle{definition}

\newtheorem{remark}       [theorem]  {Remark}

\newtheorem{application}      [theorem]  {Application}

\def \r{\mbox{${\mathbb R}$}}
\def \s{\mbox{${\mathbb S}$}}

\def \h{\mbox{${\mathbb H}$}}

\def \onabla{\mbox{$\overline{\nabla}$}}
\def \og{\mbox{$\overline{g}$}}
\def \ophi{\mbox{$\overline{\phi}$}}

\DeclareMathOperator{\trace}{trace}

\DeclareMathOperator{\grad}{grad}

\DeclareMathOperator{\Div}{div}

\DeclareMathOperator{\riem}{Riem}

\DeclareMathOperator{\vol}{Vol}

\DeclareMathOperator{\rank}{rank}

\DeclareMathOperator{\sym}{sym}
\DeclareMathOperator{\Hess}{Hess}

\numberwithin{equation}{section}

\begin{document}

\title{The stress-energy tensor for biharmonic maps }

\author{E.~Loubeau}
\author{S.~Montaldo}
\author{C.~Oniciuc}

\address{D{\'e}partement de Math{\'e}matiques \\
Laboratoire CNRS UMR 6205 \\
Universit{\'e} de Bretagne Occidentale \\
6, avenue Victor Le Gorgeu \\
CS 93837, 29238 Brest Cedex 3, France}
\email{loubeau@univ-brest.fr}

\address{Universit\`a degli Studi di Cagliari\\
Dipartimento di Matematica e Informatica\\
Via Ospedale 72\\
09124 Cagliari, Italia}
\email{montaldo@unica.it}

\address{Faculty of Mathematics\\ ``Al.I. Cuza'' University of Iasi\\
Bd. Carol I no. 11 \\
700506 Iasi, Romania}
\email{oniciucc@uaic.ro}

\begin{abstract}
Using Hilbert's criterion, we consider the stress-energy tensor
associated to the bienergy functional. We show that it derives
from a variational problem on metrics and exhibit the peculiarity
of dimension four. First, we use this tensor  to construct new examples
of biharmonic maps, then  classify maps with vanishing or
parallel stress-energy tensor and Riemannian immersions whose
stress-energy tensor is proportional to the metric.
\end{abstract}

\date{}



\thanks{
The third author was supported by a CNR-NATO (Italy) fellowship}

\maketitle

\section{Introduction}

As described by Hilbert in~\cite{DH}, the {\it stress-energy}
tensor associated to a variational problem is a symmetric
$2$-covariant tensor $S$ conservative at critical points,
i.e. $\Div S=0$.

In the context of harmonic maps, i.e. critical points of the
energy $E(\phi)=\frac{1}{2}\int_M\vert d\phi\vert^2 \ v_g$,
 the stress-energy tensor was studied in details by
Baird and Eells in~\cite{PBJE}. Indeed, the Euler-Lagrange
equation associated to the energy is the vanishing of the tension
field $\tau(\phi)=\trace\nabla d\phi$, and the tensor
$$
S=\frac{1}{2}\vert d\phi\vert^2 g - \phi^{\ast}h
$$
satisfies $\Div S=-\langle\tau(\phi),d\phi\rangle$.

As shown by Sanini in~\cite{AS}, $S$ vanishes precisely at
critical points of the energy for variations of the domain metric,
rather than variations of the map.

In this paper we consider a natural generalization of harmonic
maps, suggested by Eells and Sampson~\cite{JEJHS}: the  {\it
bienergy} of $\phi:(M,g)\to (N,h)$ is
$$
E_2(\phi)=\frac{1}{2}\int_M\vert\tau(\phi)\vert^2 \ v_g,
$$
and a map is {\it biharmonic} if it is a critical point of $E_2$,
equivalently, if it satisfies the associated Euler-Lagrange
equation
$$
\tau_2(\phi)=-\Delta\tau(\phi) -\trace R^N(d\phi,\tau(\phi))d\phi
= 0.
$$
In \cite{GYJ2}, Jiang constructed an ad-hoc $(0,2)$-tensor
\begin{eqnarray*}
S_2(X,Y)&=&\frac{1}{2}\vert\tau(\phi)\vert^2\langle X,Y\rangle+
\langle d\phi,\nabla\tau(\phi)\rangle \langle X,Y\rangle \\
\nonumber && -\langle d\phi(X), \nabla_Y\tau(\phi)\rangle-\langle
d\phi(Y), \nabla_X\tau(\phi)\rangle,
\end{eqnarray*}
such that $\Div S_2=-\langle\tau_2(\phi),d\phi\rangle$, thus
conforming to the principle of a  stress-energy tensor for the
bienergy. In analogy with harmonic maps, one would expect that
such a tensor could be depicted by metric variations.
Theorem~\ref{eq:teorema-prima-variatie} shows that $S_2$
does indeed possess this property.

Motivated by this characterization, we study $S_2=0$ in details,
and describe, in Section~\ref{eq:S-2=0}, situations for which this
implies harmonicity of the map. On compact domains ($m\neq 4$),
this was already proved by Jiang. This specialness of dimension
four is the subject of Theorem~\ref{eq:imersie-pseudo-umbilicala}
and Proposition~\ref{eq:imersie-conforma}, where Riemannian and
conformal immersions with $S_2=0$ are classified.

A salient feature of this tensor is its usefulness in finding new
biharmonic maps. In fact, the search for maps with vanishing
bitension field is replaced with one for divergence-free
stress-energy tensor, as implemented in
Proposition~\ref{eq:aplicatie-S-2}
and~\ref{eq:biarmonicitate-identitate}.

Since $S_2=0$ is a particularly strong condition, we study, in
Section ~\ref{eq:S-2-paralel}, maps with parallel stress-energy
tensor. This clearly is stronger than being simply divergence
free. In Theorem ~\ref{eq:S-2-paralel-1}, we classify
hypersurfaces with parallel stress-energy tensor and, in
Proposition~\ref{eq:S-2-paralel-2} characterize pseudo-umbilical
Riemannian immersions with parallel stress-energy tensor.

Finally, we point out that, if $M$ is compact, the bienergy is
homogeneous of degree zero with respect to the metric $g$ if and
only if $m=4$. Therefore, if $m\neq 4$, changing homothetically
the domain metric, we can render the bienergy arbitrarily large or
small. This leads to considering volume-preserving variations of
the metric, a natural problem in geometry as for the
characterization of Einstein metrics \cite{MB}. The last section
classifies Riemannian immersions  which are critical points of
$E_2$ with respect to isovolumetric variations.

See~\cite{SMCO} for an account of biharmonic maps
and~\cite{BMBib} for an up-to-date bibliography.

\noindent {\it Conventions}. We work in the $C^{\infty}$ category,
i.e. manifolds, metrics, connections and maps are
smooth, and $(M^m,g)$  denotes a connected manifold of
dimension $m$, without boundary, endowed with a Riemannian metric
$g$. By an abuse of notation, $\langle,\rangle$ indicates
inner products from different vector bundles. The Riemann
curvature operator is defined by
$R(X,Y)Z=[\nabla_{X},\nabla_{Y}]Z-\nabla_{[X,Y]}Z$. The notations
$\sharp$ and $\flat$ are for the standard musical isomorphisms.
Compactness is not assumed unless explicitly stated.

\section{The Euler-Lagrange equation and applications}

For smooth maps $\phi:(M,g)\to (N,h)$ between Riemannian
manifolds, $M$ compact and orientable,  consider the
bienergy functional $E_2$:
$$
E_2:C^{\infty}(M,N)\to \r, \quad E_2(\phi)=\int_M
\vert\tau(\phi)\vert^2 \ v_g,
$$
and, as mentioned in the introduction, a map is biharmonic if it
is a critical point of $E_2$, that is, for any variation
$\{\phi_t\}$ of $\phi$,
$\frac{d}{dt}\big{\vert}_{t=0}E_2(\phi_t)=0$.

Opting for a different angle of attack, one can vary the
metric instead of the map, more precisely, given $\phi:M\to (N,h)$,
consider the functional
$$
F:G\to \r, \quad F(g)=E_2(\phi),
$$
where $G$ is the set of Riemannian metrics on $M$. As $G$ is an
infinite dimensional manifold (\cite{AB}), it admits a tangent space
at $g$, the set of symmetric $(0,2)$-tensors on $M$, i.e.
$$
T_gG=C(\odot^2T^{\ast}M).
$$

\noindent For a curve $t\to g_t$ in $G$ with $g_0=g$, denote by
$$
\omega=\frac{d}{dt}\big{\vert}_{t=0}\{g_t\}=\delta(g_t)\in T_gG
$$
the corresponding variational tensor field which, in local
coordinates, can be written
$$
\omega=\frac{\partial g_{ij}}{\partial
t}(x,0)dx^idx^j=\omega_{ij}dx^idx^j,
$$
where $g_t=g_{ij}(x,t)dx^idx^j$, and write
$\delta = \frac{d}{dt}\big{\vert}_{t=0}$ for the first variation.
\newline For a one-parameter variation $\{g_t\}$ of $g$ we have
$$
F(g_t)=\frac{1}{2}\int_M \vert\tau_t(\phi)\vert^2 \ v_{g_t}.
$$
We now compute $\delta(F(g_t))$. Differentiating $F(g_t)$
leads to:
\begin{equation}
\label{eq:derivata1} \delta(F(g_t))=\frac{1}{2}\int_M
\delta(\vert\tau_t(\phi)\vert^2) \ v_g + \frac{1}{2}\int_M
\vert\tau(\phi)\vert^2 \ \delta(v_{g_t}).
\end{equation}
The calculation of the first term breaks down in two lemmas.

\begin{lemma}
\label{eq:lema1} The vector field $\xi=(\Div
\omega)^{\sharp}-\frac{1}{2}\grad(\trace\omega)$ satisfies:
$$
\delta(\vert\tau_t(\phi)\vert^2)=-2\langle \tau(\phi).\nabla
d\phi,\omega\rangle-2\langle\tau(\phi),d\phi(\xi)\rangle,
$$
where $\tau(\phi).\nabla d\phi\in C(\odot^2T^{\ast}M)$ is intended as
$$
\big(\tau(\phi).\nabla d\phi\big)(X,Y)=\langle\tau(\phi),\nabla
d\phi(X,Y)\rangle.
$$
\end{lemma}

\begin{proof}
In local coordinates $\{(U;x^i)\}_{i=1}^m$ on $M$ and
$\{(V;y^{\alpha})\}_{\alpha=1}^n$ on $N$:
\begin{equation}
\label{eq:derivata-tau-1} \delta(\vert\tau_t(\phi)\vert^2)=
\delta(h_{\alpha\beta}(\phi)\tau^{\alpha}_t\tau^{\beta}_t)=
2h_{\alpha\beta}(\delta(\tau^{\alpha}_t))\tau^{\beta}.
\end{equation}

\noindent Now
\begin{eqnarray*}
\delta(\tau^{\alpha}_t)&=&
\delta\Big(g^{ij}(x,t)\big(\frac{\partial^2\phi^{\alpha}}{\partial
x^i\partial x^j}-\Gamma^k_{ij}(x,t)\phi^{\alpha}_k
+^N\Gamma^{\alpha}_{\beta\sigma}\phi^{\beta}_i\phi^{\sigma}_j
\big)
\Big) \\
&=&(\delta g^{ij})(\nabla
d\phi)^{\alpha}_{ij}-g^{ij}(\delta\Gamma^k_{ij})\phi^{\alpha}_k.
\end{eqnarray*}

\noindent Since $g^{il}(x,t)g_{lk}(x,t)=\delta^i_k$, 
$
\delta g^{ij}=-g^{ia}g^{jb}(\delta
g_{ab})=-g^{ia}g^{jb}\omega_{ab},
$
so
\begin{equation}
\label{eq:derivata-tau-2}
\delta(\tau^{\alpha}_t)=-g^{ia}g^{jb}\omega_{ab} (\nabla
d\phi)^{\alpha}_{ij}-g^{ij}(\delta\Gamma^k_{ij})\phi^{\alpha}_k.
\end{equation}

\noindent Now we compute the term
$g^{ij}(\delta\Gamma^k_{ij})\phi^{\alpha}_k$. From
$$
\Gamma^k_{ij}(x,t)=\frac{1}{2}g^{kl}(x,t)\big(\frac{\partial
g_{li}}{\partial x^j}(x,t)+\frac{\partial g_{lj}}{\partial
x^i}(x,t)-\frac{\partial g_{ij}}{\partial x^l}(x,t)\big)
$$
we have:
\begin{eqnarray}
\label{eq:derivata-gamma-1} \nonumber \delta\Gamma^k_{ij}&=&
\frac{1}{2}(\delta g^{kl})\big(\frac{\partial g_{li}}{\partial
x^j}+\frac{\partial g_{lj}}{\partial x^i}-\frac{\partial
g_{ij}}{\partial x^l}\big) \\
\nonumber &&+ \frac{1}{2} g^{kl}\big(\frac{\partial^2
g_{li}}{\partial x^j
\partial t}(x,0)+\frac{\partial^2 g_{lj}}{\partial x^i\partial
t}(x,0)-\frac{\partial^2 g_{ij}}{\partial x^l\partial t}(x,0)\big) \\
&=&-\frac{1}{2}g^{ka}g^{lb}\omega_{ab}\big(\frac{\partial
g_{li}}{\partial x^j}+\frac{\partial g_{lj}}{\partial
x^i}-\frac{\partial g_{ij}}{\partial x^l}\big) + \frac{1}{2}
g^{kl}\big(\frac{\partial \omega_{li}}{\partial
x^j}+\frac{\partial \omega_{lj}}{\partial
x^i}-\frac{\partial \omega_{ij}}{\partial x^l}\big) \\
\nonumber &=&-g^{ka}\Gamma^b_{ij}\omega_{ab}+\frac{1}{2}
g^{kl}\big(\frac{\partial \omega_{li}}{\partial
x^j}+\frac{\partial \omega_{lj}}{\partial
x^i}-\frac{\partial \omega_{ij}}{\partial x^l}\big). \\
\nonumber
\end{eqnarray}
But
$$
\nabla_j\omega_{li}=\frac{\partial\omega_{li}}{\partial x^j}
-\Gamma^h_{jl}\omega_{hi}-\Gamma^h_{ji}\omega_{lh}
$$
implies
$$
\frac{\partial \omega_{li}}{\partial x^j}+\frac{\partial
\omega_{lj}}{\partial x^i}-\frac{\partial \omega_{ij}}{\partial
x^l}=\nabla_j\omega_{li}+\nabla_i\omega_{lj}-\nabla_l\omega_{ij}
+2\Gamma^h_{ij}\omega_{hl}.
$$
So replacing in Equation~\eqref{eq:derivata-gamma-1} we
obtain:
\begin{equation}
\label{eq:derivata-gamma-2} \delta\Gamma^k_{ij}=
\frac{1}{2}g^{kl}(\nabla_j\omega_{li}+\nabla_i\omega_{lj}-\nabla_l\omega_{ij})
\end{equation}
and
\begin{equation}
\label{eq:derivata-gamma-3}
g^{ij}(\delta\Gamma^k_{ij})\phi^{\alpha}_k=
g^{kl}(\nabla_j\omega^j_l-\frac{1}{2}\nabla_l\trace\omega)\phi^{\alpha}_k=
\xi^k\phi^{\alpha}_k.
\end{equation}
From ~\eqref{eq:derivata-tau-1}, ~\eqref{eq:derivata-tau-2} and
~\eqref{eq:derivata-gamma-3}, the lemma follows.
\end{proof}

\begin{lemma}
\label{eq:lema2} Consider the one-form $d\phi.\tau(\phi)\in
\Lambda^1(M)$ defined by $d\phi.\tau(\phi)(X)=\langle
d\phi(X),\tau(\phi)\rangle$, and
$\sym\big(\nabla(d\phi.\tau(\phi))\big)$ the symmetric part of
$\nabla(d\phi.\tau(\phi))$, then:
$$
\int_M \langle \tau(\phi),d\phi(\xi)\rangle \ v_g= \int_M
\langle-\sym\big(\nabla(d\phi.\tau(\phi))\big)+\frac{1}{2}\big(\Div\big(d\phi.\tau(\phi)
\big)^{\sharp}\big)g,\omega\rangle \ v_g .
$$
\end{lemma}

\begin{proof}
First observe that:
\begin{equation}
\label{eq:lema2-1} \langle\tau(\phi),d\phi(Z)\rangle=\langle
d\phi.\tau(\phi),Z^{\flat}\rangle, \quad \forall Z\in C(TM).
\end{equation}
By the definition of $\xi$
\begin{eqnarray}
\label{eq:lema2-2} \int_M \langle\tau(\phi),d\phi(\xi)\rangle \
v_g&=& \int_M
\langle\tau(\phi),d\phi\big((\Div\omega)^{\sharp}\big)\rangle \
v_g \\
\nonumber && -\frac{1}{2}\int_M \langle\tau(\phi),d\phi\big(\grad
(\trace \omega)\big)\rangle \  v_g
\end{eqnarray}
and, by~\eqref{eq:lema2-1}, the first term on the right-hand
side of ~\eqref{eq:lema2-2} becomes
$$
\int_M
\langle\tau(\phi),d\phi\big((\Div\omega)^{\sharp}\big)\rangle \
v_g = \int_M \langle d\phi.\tau(\phi), \Div\omega\rangle \ v_g.
$$
On the other hand, if $\theta\in \Lambda^1(M)$, $\sigma\in
C(\odot^2T^{\ast}M)$, and
$C(\theta,\sigma)=(\theta^i\sigma_{ij})dx^j=(\theta_i\sigma^i_j)dx^j$
denotes their contraction, we have:
\begin{equation}
\label{eq:formula-divergenta-1}
\langle\theta,\Div\sigma\rangle=\Div(C(\theta,\sigma)^{\sharp})-
\langle\sym(\nabla\theta),\sigma\rangle.
\end{equation}
Applying ~\eqref{eq:formula-divergenta-1} to
$\theta=d\phi.\tau(\phi)$ and $\sigma=\omega$, yields
\begin{equation}
\label{eq:lema2-3} \int_M
\langle\tau(\phi),d\phi\big((\Div\omega)^{\sharp}\big)\rangle \
v_g = -\int_M\langle \sym\big(\nabla(d\phi.\tau(\phi))\big),\omega
\rangle \ v_g.
\end{equation}

\noindent The second term on the right-hand side of
~\eqref{eq:lema2-2} can then be written:
\begin{eqnarray}
\nonumber \label{eq:lema2-4} \int_M \langle
\tau(\phi),d\phi\big(\grad (\trace \omega)\big)\rangle \
v_g&=&\int_M \langle
d\phi.\tau(\phi), d(\trace \omega)\rangle \ v_g \\
\nonumber &=& \int_M \langle \big(d\phi.\tau(\phi)\big)^{\sharp},
\grad(\trace \omega)\rangle \
v_g \\
&=& -\int_M (\trace
\omega)
\Div\big(d\phi.\tau(\phi)\big)^{\sharp} \ v_g \\
\nonumber &=& -\int_M \langle
\big(\Div\big(d\phi.\tau(\phi)\big)^{\sharp}\big)g,\omega\rangle \
v_g.
\end{eqnarray}
The lemma follows from ~\eqref{eq:lema2-2}, ~\eqref{eq:lema2-3}
and ~\eqref{eq:lema2-4}.
\end{proof}

This preparation is the key to:

\begin{theorem}
\label{eq:teorema-prima-variatie}
Let $\phi:(M,g)\to (N,h)$ be a
smooth map, $M$ compact and orientable, and $\{g_t\}$  a
one-parameter variation of $g$ through Riemannian metrics. Then
$$
\delta(F(g_t))=-\frac{1}{2}\int_M \langle S_2,\omega\rangle \ v_g,
$$
where $S_2\in C(\odot^2T^{\ast}M)$ is given by:
\begin{eqnarray}
\label{eq:formula-S2}
S_2(X,Y)&=&\frac{1}{2}\vert\tau(\phi)\vert^2\langle X,Y\rangle+
\langle d\phi,\nabla\tau(\phi)\rangle \langle X,Y\rangle \\
\nonumber && -\langle d\phi(X), \nabla_Y\tau(\phi)\rangle-\langle
d\phi(Y), \nabla_X\tau(\phi)\rangle.
\end{eqnarray}
\end{theorem}

\begin{proof}
Recall that 
$
\delta(v_{g_t})=\langle\frac{1}{2}g,\omega\rangle v_g  
$
(see, for example, \cite{AS,PBJCW}).
Then, by Lemma
~\ref{eq:lema1} and ~\ref{eq:lema2}, we can rewrite
~\eqref{eq:derivata1}:
\begin{eqnarray*}
\delta(F(g_t))&=&\int_M
\Big\langle\big[\frac{1}{4}\vert\tau(\phi)\vert^2g
-\frac{1}{2}\big(\Div\big(d\phi.\tau(\phi)\big)^{\sharp}\big)g \\
&& \ \ \ \ \ \
+\sym\big(\nabla(d\phi.\tau(\phi))\big)-\tau(\phi).\nabla
d\phi\big],\omega\Big\rangle  \ v_g.
\end{eqnarray*}
The formula (cf.~\cite{JELL})
\begin{equation}
\label{eq:formula-divergenta-2} \Div
\big(d\phi.\tau(\phi)\big)^{\sharp}=\vert
\tau(\phi)\vert^2+\langle d\phi,\nabla\tau(\phi)\rangle
\end{equation}
and the expression of $\sym\big(\nabla(d\phi.\tau(\phi))\big)$:
\begin{eqnarray*}
\sym\big(\nabla(d\phi.\tau(\phi))\big)(X,Y)&=&
\frac{1}{2}\Big(\big(\nabla(d\phi.\tau(\phi))\big)(X,Y)
+\big(\nabla(d\phi.\tau(\phi))\big)(Y,X)\Big) \\
&=&\langle\nabla d\phi(X,Y),\tau(\phi)\rangle \\
&& + \frac{1}{2}\big(\langle d\phi(X),\nabla_Y\tau(\phi)\rangle
+\langle d\phi(Y),\nabla_X\tau(\phi)\rangle\big),
\end{eqnarray*}
end the proof.
\end{proof}

As mentioned in the introduction, $S_2$ was known by Jiang
in~\cite{GYJ2}, where he proved the following

\begin{theorem}
[\cite{GYJ2}] \label{eq:teorema-divergenta} For any map $\phi:(M,g)\to
(N,h)$:
$$
\Div S_2(Y)=-\langle\tau_2(\phi),d\phi(Y)\rangle, \quad \forall \
Y\in C(TM).
$$
\end{theorem}

\begin{proof}
We give a proof for the sake of completeness. Write
$S_2=T_1+T_2$, where $T_1,T_2\in C(\odot^2T^{\ast}M)$
are defined by
$$\begin{array}{l}
T_1(X,Y)=\frac{1}{2}\vert\tau(\phi)\vert^2\langle X,Y\rangle+
\langle d\phi,\nabla\tau(\phi)\rangle \langle X,Y\rangle
\\
\\
T_2(X,Y)=-\langle d\phi(X), \nabla_Y\tau(\phi)\rangle-\langle
d\phi(Y), \nabla_X\tau(\phi)\rangle.
\end{array}
$$
Let $p\in M$ and $\{X_i\}_{i=1}^m$ a
geodesic frame centered on $p$. Writing $Y=Y^iX_i$, at $p$, we have:
\begin{eqnarray}
\label{eq:divergenta1} \nonumber \Div
T_1(Y)&=&\sum_i(\nabla_{X_i}T_1)(X_i,Y)=
\sum_i\big(X_i(T_1(X_i,Y))-T_1(X_i,\nabla_{X_i}Y)\big) \\
\nonumber &=&\sum_i\Big(
X_i\big(\frac{1}{2}\vert\tau(\phi)\vert^2Y^i+\sum_j\langle
d\phi(X_j),\nabla_{X_j}\tau(\phi)\rangle Y^i\big) \\
&&-\frac{1}{2}\vert\tau(\phi)\vert^2X_iY^i- \sum_j\langle
d\phi(X_j),\nabla_{X_j}\tau(\phi)\rangle(X_iY^i) \Big) \\
\nonumber
&=&\nonumber\langle\nabla_Y\tau(\phi),\tau(\phi)\rangle
+\sum_i\langle\nabla d\phi(Y,X_i),\nabla_{X_i}\tau(\phi)\rangle \\
\nonumber && +\sum_i\langle
d\phi(X_i),\nabla_Y\nabla_{X_i}\tau(\phi)\rangle \\
&=&\nonumber\langle\nabla_Y\tau(\phi),\tau(\phi)\rangle+\trace\langle\nabla
d\phi(Y,\cdot),\nabla_{\cdot}\tau(\phi)\rangle \\
\nonumber &&+\trace\langle
d\phi(\cdot),\nabla^2\tau(\phi)(Y,\cdot)\rangle,
\end{eqnarray}
whilst
\begin{eqnarray}
\label{eq:divergenta2} \nonumber \Div
T_2(Y)&=&\sum_i\big(X_i(T_2(X_i,Y))-T_2(X_i,\nabla_{X_i}Y)\big) \\
\nonumber &=&-\langle\nabla_Y\tau(\phi),\tau(\phi)\rangle -
\sum_i\langle\nabla d\phi(Y,X_i),\nabla_{X_i}\tau(\phi)\rangle \\
&&-\sum_i\langle
d\phi(X_i),\nabla_{X_i}\nabla_Y\tau(\phi)
-\nabla_{\nabla_{X_i}Y}\tau(\phi)\rangle
\\
\nonumber &&+\langle d\phi(Y),\Delta\tau(\phi)\rangle \\
\nonumber
&=&-\langle\nabla_Y\tau(\phi),\tau(\phi)\rangle-\trace\langle\nabla
d\phi(Y,\cdot),\nabla_{\cdot}\tau(\phi)\rangle \\
\nonumber &&-\trace\langle
d\phi(\cdot),\nabla^2\tau(\phi)(\cdot,Y)\rangle +\langle
d\phi(Y),\Delta\tau(\phi)\rangle.
\end{eqnarray}
Summing~\eqref{eq:divergenta1} and~\eqref{eq:divergenta2} gives:
\begin{eqnarray*}
\Div S_2(Y)&=&\langle d\phi(Y),\Delta\tau(\phi)\rangle+
\sum_i\langle
d\phi(X_i),R(Y,X_i)\tau(\phi)\rangle \\
&=&
-\langle\tau_2(\phi),d\phi(Y)\rangle.
\end{eqnarray*}
\end{proof}

Theorem~\ref{eq:teorema-divergenta} links $S_2$ with the bitension
field and immediately leads to the following characterizations.

\begin{corollary}
\label{eq:corolar-1} Let $\phi:(M,g)\to (N,h)$ be a smooth map.
\begin{itemize}
\item[a)] If $\phi$ is a Riemannian immersion, then $\Div S_2=0$
if and only if the tangent part of $\tau_2(\phi)$ vanishes.
\item[b)] If $\phi$ is a submersion (not necessarily Riemannian),
then $\Div S_2=0$ if and only if $\phi$ is biharmonic.
\end{itemize}
\end{corollary}

These results are applied to construct proper (i.e. non-harmonic) biharmonic maps.

\begin{proposition}
\label{eq:aplicatie-S-2} Let $\phi:(M,g)\to (N,h)$ be a
submersion, $M$ non-compact. Assume that $\tau(\phi)$ is basic,
i.e. $\tau(\phi)=\zeta\circ\phi$, where $\zeta\in C(TN)$. If
$\zeta$ is Killing and $\vert\zeta\vert^2=c^2\neq 0$ is constant,
then $\phi$ is proper biharmonic.
\end{proposition}

\begin{proof}
We prove that $\Div S_2=0$. The tensor $S_2$
takes, in this case, the expression:
$$
S_2(X,Y)=\big\{\frac{c^2}{2} + \langle
d\phi,\nabla\tau(\phi)\rangle\big\} \langle X,Y\rangle -\langle
d\phi(X), \nabla_Y\tau(\phi)\rangle-\langle d\phi(Y),
\nabla_X\tau(\phi)\rangle.
$$

\noindent Now, let $p$ be a point in $M$ and
$\{X_i\}_{i=1}^m$ an orthonormal basis of $T_pM$ such that
$\{X_{\alpha}\}_{\alpha=1}^n$ belongs to
$T^H_pM=(T^V_pM)^{\perp}$ and $\{X_a\}_{a=n+1}^m$ to $T^V_pM=\ker d\phi(p)$.
As $\zeta$ is Killing:
\begin{eqnarray*}
\langle d\phi,\nabla\tau(\phi)\rangle (p)&=&\sum_{\alpha}\langle
d\phi_p(X_{\alpha}),\nabla_{X_{\alpha}}\tau(\phi)\rangle+
\sum_a\langle d\phi_p(X_a),\nabla_{X_a}\tau(\phi)\rangle \\
&=&\sum_{\alpha}\langle
d\phi_p(X_{\alpha}),\nabla^N_{d\phi_p(X_{\alpha})}\zeta\rangle \\
&=&0
\end{eqnarray*}
and, for any $X,Y\in T_pM$,
\begin{eqnarray*}
S_2(p)(X,Y)&=&\frac{c^2}{2}\langle X,Y\rangle-\big(\langle
d\phi_p(X),\nabla^N_{d\phi_p(Y)}\zeta\rangle + \langle
d\phi_p(Y),\nabla^N_{d\phi_p(X)}\zeta\rangle\big) \\
&=&\frac{c^2}{2}\langle X,Y\rangle.
\end{eqnarray*}
Thus $S_2=\frac{c^2}{2}g$ is divergence free.

If $M$ is compact, since $\langle d\phi,\nabla\tau(\phi)\rangle=0$,
integrating~\eqref{eq:formula-divergenta-2}, yields
$\tau(\phi)=0$.
\end{proof}

\noindent We apply Proposition~\ref{eq:aplicatie-S-2} to three
situations. Recall that a function is affine if its restriction to
any geodesic is affine with respect to the real parameter, or,
equivalently, if its gradient is parallel (\cite{TS}).

\begin{application}
\label{eq:primul-exemplu} Let $(M,g)$ and $(N,h)$ be Riemannian
manifolds and denote by $M\times_{f^2}N$ their warped product with
respect to a positive function on $M$, then the projection
$\pi:M\times_{f^2}N\to M$ is a Riemannian submersion with
$\tau(\pi)=n\grad(\ln f)\circ\pi$. When $\ln f$ is affine,
$\grad(\ln f)$ is Killing of constant norm (\cite{ABSMCO}), hence
$\pi$ is biharmonic.
\end{application}

\begin{application}
For a vector field $\zeta$ on $M$, $TM$ can be endowed
with a Sasaki-type metric such that the projection
$\pi:TM\to M$ is a Riemannian submersion and
$\tau(\pi)=-(m+1)\zeta\circ\pi$ (\cite{CO}).
Choosing $\zeta$ Killing of constant norm makes $\pi$ biharmonic.
\end{application}

In Application~\ref{eq:primul-exemplu}, the fibres of $\pi$  do
not have the same volume, otherwise the next proposition
would force the map to be harmonic.

\begin{proposition}
Let $\phi:(M,g)\to (N,h)$ be a Riemannian submersion with basic
tension field, $M$ compact and orientable. If $\vol(\phi^{-1}(q))$
is constant then $\phi$ is harmonic.
\end{proposition}

\begin{proof}
Write $\tau(\phi)=\zeta\circ\phi$, $\zeta\in C(TN)$ then, by a
simple calculation, $\langle
d\phi,\nabla\tau(\phi)\rangle=(\Div\zeta)\circ\phi$. If
$\vol(\phi^{-1}(q))=c$ for all $q\in N$, the coarea formula
(\cite{MBor,YDBVAZ}), would impose:
\begin{eqnarray*}
\int_M\langle d\phi,\nabla\tau(\phi)\rangle \ v_g&=&
\int_M(\Div\zeta)\circ\phi \ v_g \\
&=&\int_N \Big(\int_{\phi^{-1}(q)}(\Div\zeta)(q) \ v_{g_q}\Big) \
v_h  =c\int_N\Div\zeta \ v_h=0,
\end{eqnarray*}
where $g_q$ is the induced metric on $\phi^{-1}(q)\subset M$. On
the other hand,
$$
\int_M\langle d\phi,\nabla\tau(\phi)\rangle \ v_g
=-\int_M\vert\tau(\phi)\vert^2 \
v_g=-\int_M\vert\zeta\vert^2\circ\phi \
v_g=-c\int_N\vert\zeta\vert^2 \ v_h .
$$
\end{proof}

\begin{application}
A submersion $\phi:(M,g)\to (N,h)$ is horizontally conformal
(\cite{PBJCW}) if there exists a positive function $\lambda\in
C(M)$, the dilation, such that
$$
h(d\phi_p(X),d\phi_p(Y))=\lambda^2g_p(X,Y),
\quad \forall X,Y\in (\ker d\phi_p)^{\perp}.
$$

\noindent A partial biharmonic analogue of the Baird-Eells Theorem, which
states that, for horizontally homothetic submersion harmonicity is
equivalent to minimality of the fibres, would be:

\begin{proposition}
\label{eq:submersie-orizontal-conforma} Let $\phi:(M,g)\to (N,h)$
be a submersion with dilation $\lambda=\tilde{\lambda}\circ\phi$
and minimal fibers. If $\tilde{\lambda}^2$ is non-constant and
affine, and $n\neq 2$, then $\phi$ is proper biharmonic.
\end{proposition}

\begin{proof}
By a straightforward computation we get
$\tau(\phi)=\frac{2-n}{2}(\grad \tilde{\lambda}^2)\circ\phi$ and,
since $\tilde{\lambda}^2$ is affine, then $\grad\tilde{\lambda}^2$
is Killing of constant norm.
\end{proof}

We point out that a horizontally conformal submersion $\phi:(M,g)\to(N,h)$, with
constant dilation along the fibres, is the composition between a
Riemannian submersion $\varphi:(M,g)\to(N,\overline{h})$ and the identity map 
${\bf 1}:(N,\overline{h})\to(N,h)$, where $\overline{h}$ is conformal to $h$.
\end{application}

\begin{proposition}
\label{eq:biarmonicitate-identitate} Let ${\bf 1}:(M^m,g)\to
(M^m,\og=e^{2\rho}g)$ be the identity map, $m\neq 2$ and $\rho\in
C^{\infty}(M)$.
\begin{itemize}
\item[a)] If $\rho=\ln\sqrt{f}$, $f$ a nonconstant affine positive
function on $(M,\og)$, then ${\bf 1}$ is proper biharmonic.
\item[b)] If $\rho$ is a nonconstant affine function on $(M,\og)$,
then ${\bf 1}$ is proper biharmonic if and only if $m=4$.
\end{itemize}
\end{proposition}

\begin{proof} a)
To prove $\Div S_2=0$, we compute:
$$
\tau({\bf 1})=(2-m)\grad_g\rho=(2-m)e^{2\rho}\grad_{\og}\rho=
\frac{2-m}{2}\grad_{\og}(e^{2\rho}).
$$
By hypothesis, $e^{2\rho}=f$ is affine on $(M,\og)$,
thus $\grad_{\og}(e^{2\rho})$ is parallel and
$\onabla\tau({\bf 1})=0$, $\onabla$ being the connection of $\og$.
Therefore $S_2=\frac{1}{2}c^2g$,
with $c^2=\vert\tau({\bf 1})\vert^2_{\og}$ constant.

\noindent We note that this result can be also obtained by
Proposition~\ref{eq:submersie-orizontal-conforma}.

\noindent b) From $\onabla\grad_{\og}\rho=0$ we infer:
$$
\vert\tau({\bf 1})\vert^2_{\og}=(2-m)^2c^2e^{4\rho}, \quad
\onabla_{X}\tau({\bf 1})=2(2-m)e^{2\rho}(X\rho)\grad_{\og}\rho,
$$
where $c^2=\vert\grad_{\og}\rho\vert^2_{\og}$ is constant. Put
$S_2=T_1-T_2$, with
$$
T_1=\big(\frac{1}{2}\vert\tau({\bf 1})\vert^2+\langle d{\bf
1},\nabla\tau({\bf 1})\rangle\big)g
=\frac{(2-m)(6-m)}{2}c^2e^{4\rho}g
$$
and
$$
T_2(X,Y)=\langle d{\bf 1}(X),\nabla_Y\tau({\bf 1})\rangle+ \langle
d{\bf 1}(Y),\nabla_X\tau({\bf
1})\rangle=4(2-m)e^{2\rho}(X\rho)(Y\rho).
$$
Then:
\begin{eqnarray}
\label{eq:prima-divergenta} \Div
T_1(Z)&=&\frac{(2-m)(6-m)}{2}c^2g(\grad e^{4\rho},Z)
\\
\nonumber
&=&2(2-m)(6-m)c^2e^{4\rho}(Z\rho),
\end{eqnarray}
and
\begin{eqnarray}
\label{eq:a-doua-divergenta-1} \Div
T_2(Z)&=&8(2-m)c^2e^{4\rho}(Z\rho)-4(2-m)e^{2\rho}(\Delta\rho)(Z\rho)
\\
\nonumber &&+4(2-m)e^{2\rho}\Hess\rho(\grad_g\rho,Z).
\end{eqnarray}
If we look at $\rho$ as $\overline{\rho}\circ {\bf 1}$, where
$\overline{\rho}:(M,\og)\to \r$, and apply the chain rule, we get
\begin{equation}
\label{eq:a-doua-divergenta-2} -\Delta\rho=(2-m)c^2e^{2\rho},
\quad \Hess\rho(\grad_g\rho,Z)= c^2e^{2\rho}(Z\rho).
\end{equation}
Replace \eqref{eq:a-doua-divergenta-2} in
~\eqref{eq:a-doua-divergenta-1} to obtain:
\begin{equation}
\label{eq:a-doua-divergenta} \Div
T_2(Z)=4(2-m)(5-m)c^2e^{4\rho}(Z\rho),
\end{equation}
and from ~\eqref{eq:prima-divergenta} and
~\eqref{eq:a-doua-divergenta}
$$
\Div S_2(Z)=2(2-m)(m-4)c^2e^{4\rho}(Z\rho).
$$
We conclude since $\rho$ is not constant and $m\neq 2$.
\end{proof}

\begin{remark}
The function $\rho$ in Proposition
~\ref{eq:biarmonicitate-identitate} is isoparametric (though not
affine) on $(M,g)$. On Einstein manifolds, interesting links
between biharmonicity of the identity (modulo a conformal change
of the metric) and isoparametricity of the conformal factor have
been discovered in~\cite{PBDK,ABA}.
\end{remark}

We close this section with another application of $\Div S_2=0$.

\noindent Let $\sigma\in C(\odot^2 T^{\ast}M)$,
$\sigma=\sigma_{ij}dx^idx^j$, and $X\in C(TM)$,
$X=X^i\frac{\partial}{\partial x^i}$. Consider the contraction
$C(X,\sigma)=X^jg^{ih}\sigma_{hj}\frac{\partial}{\partial
x^i}=X^j\sigma^i_j\frac{\partial}{\partial
x^i}=(\sigma(X,\cdot))^{\sharp}$. As
$$
\Div C(X,\sigma)=\Div \sigma(X)+\frac{1}{2}\langle
\sigma,L_Xg\rangle
$$
we deduce:

\begin{proposition}
If $X$ is Killing and $\phi:(M,g)\to
(N,h)$ biharmonic map, then $C(X,S_2)$ is divergence free.
\end{proposition}

\section{Vanishing of the biharmonic stress-energy tensor}
\label{eq:S-2=0}

Clearly, from ~\eqref{eq:formula-S2}, harmonic implies $S_2=0$,
so it is only natural to study the converse.
Note that $F(g)$ is nonnegative and zero if and only if
$\phi$ is harmonic. Thus our quest is of critical points ($S_2$=0)
which are minima.

Before embarking on this problem, observe that $S_2=0$ does not, in
general, imply harmonicity, as illustrated by the non-geodesic cubic curve
$\gamma(t)=t^3a$, $a\in \r^n$. Yet, if we impose arc-length
parametrization, we have:

\begin{proposition}
Let $\gamma:I\subset\r\to (N,h)$ be a curve parametrized by arc-length,
assume $S_2=0$, then $\gamma$ is geodesic.
\end{proposition}

\begin{proof}
A direct computation shows:
\begin{eqnarray*}
0&=&S_2(\frac{\partial}{\partial t},\frac{\partial}{\partial
t})=\frac{1}{2}\vert\tau(\gamma)\vert^2-\langle
d\gamma(\frac{\partial}{\partial t}),
\nabla_{\frac{\partial}{\partial t}}\tau(\gamma)\rangle
\\
&=&\frac{3}{2}\vert\tau(\gamma)\vert^2.
\end{eqnarray*}
\end{proof}

When the domain is a surface, $S_2=0$ is indeed very
strong.

\begin{proposition}
\label{eq:cazul-suprafetelor} Let $\phi:(M^2,g)\to (N,h)$ be a map
from a surface, then $S_2=0$ implies $\phi$ harmonic.
\end{proposition}

\begin{proof} The trace of $S_2$ gives:
\begin{eqnarray*}
0&=&\trace S_2=\vert\tau(\phi)\vert^2+2\langle
d\tau(\phi),d\phi\rangle -2\langle d\tau(\phi),d\phi\rangle \\
&=&\vert\tau(\phi)\vert^2.
\end{eqnarray*}
\end{proof}

For the sequel, we first need a reformulation of $S_2=0$:

\begin{proposition}
Let $\phi:(M,g)\to (N,h)$, $m\neq 2$, then $S_2=0$ if and only if
\begin{equation}
\label{eq:formula-echivalenta}
\frac{1}{m-2}\vert\tau(\phi)\vert^2\langle
X,Y\rangle+\langle\nabla_X\tau(\phi),d\phi(Y)\rangle+
\langle\nabla_Y\tau(\phi),d\phi(X)\rangle=0,
\end{equation}
$\forall X,Y\in C(TM)$.
\end{proposition}

\begin{proposition}
\label{eq:rang} A map $\phi:(M^m,g)\to (N,h)$, $m>2$, with $S_2=0$
and $\rank\phi\leq m-1$ is harmonic.
\end{proposition}

\begin{proof}
Take $p\in M$, as $\rank\phi(p)\leq m-1$,
there exists a unit vector $X_p\in \ker d\phi_p$ and
for $X=Y=X_p$, ~\eqref{eq:formula-echivalenta} becomes
$\tau(\phi)(p)=0$.
\end{proof}

\begin{corollary}
Let $\phi:(M,g)\to (N,h)$ be a submersion ($m>n$), if $S_2=0$ then
$\phi$ is harmonic.
\end{corollary}

Recall the following result, originally due to Jiang

\begin{theorem}
[\cite{GYJ2}] \label{eq:teorema-Jiang} A map $\phi:(M,g)\to
(N,h)$, $m\neq 4$, with $S_2=0$, $M$ compact and orientable, is
harmonic.
\end{theorem}

\begin{proof}
Trace of $S_2$ is
\begin{equation}
\label{eq:urma-1} 0=\trace
S_2=\frac{m}{2}\vert\tau(\phi)\vert^2+m\langle
d\tau(\phi),d\phi\rangle -2\langle d\tau(\phi),d\phi\rangle,
\end{equation}
and integrating over $M$:
$$
0=\frac{4-m}{2}\int_M \vert\tau(\phi)\vert^2 \ v_g
$$
hence, as $m\neq 4$, $\phi$ is harmonic.
\end{proof}

\begin{remark}
An alternative proof is: consider the
one-parameter variation $\{g_t=(1+t)g\}$, then
$$
F(g_t)=(1+t)^{\frac{m-4}{2}}F(g).
$$
Now $S_2=0$ implies $\delta(F(g_t))=0$ hence ($m\neq 4$),
$F(g)=0$, i.e. $\phi$ is harmonic.
\end{remark}

When $M$ is not necessarily compact,
Theorem~\ref{eq:teorema-Jiang} can be rewritten for Riemannian
immersions:

\begin{proposition}
A Riemannian immersion $\phi:(M,g)\to (N,h)$ with $S_2=0$, $m\neq
4$, is minimal.
\end{proposition}

\begin{proof}
By~\eqref{eq:formula-divergenta-2},
$\langle d\tau(\phi),d\phi\rangle=-\vert\tau(\phi)\vert^2$,
and replacing in ~\eqref{eq:urma-1}  yields
$$
\frac{4-m}{2}\vert\tau(\phi)\vert^2=0,
$$
forcing $\phi$ to be minimal.
\end{proof}

The next result introduces integral conditions ensuring
that $S_2=0$ reveals harmonicity.
First we cite Yau's version of Stokes Theorem:

\begin{lemma}
[\cite{STY}] \label{eq:lemma-Yau} Let $(M^m,g)$ be a complete
Riemannian manifold and $\omega$ a smooth integrable $(m-1)$-form
defined on $M$. Then there exists a sequence of domains $B_i$ in
$M$ such that $M=\bigcup_iB_i$, $B_i\subset B_{i+1}$, and
$\lim_{i\to \infty}\int_{B_i}d\omega=0$.
\end{lemma}

\begin{theorem}
\label{eq:teorema-Yau} Let $(M,g)$, $m\neq 4$, be an orientable
complete Riemannian manifold and $\phi:(M,g)\to (N,h)$ a map with
$S_2=0$. If $\int_M\vert d\phi.\tau(\phi)\vert \ v_g < \infty$
then $\phi$ is harmonic.
\end{theorem}

\begin{proof}
For $m=2$, this is Proposition~\ref{eq:cazul-suprafetelor},
so now assume $m\neq 2$.
The trace of $S_2$ and ~\eqref{eq:formula-divergenta-2} result in:
$$
\frac{m-4}{2(m-2)}\vert\tau(\phi)\vert^2 \ v_g=(\Div X) \
v_g=d(i_{X}v_g),
$$
where $X=\big(d\phi.\tau(\phi)\big)^{\sharp}$. We now apply
Lemma~\ref{eq:lemma-Yau} to $\omega=i_Xv_g$. To compute the norm
of $\omega$, choose $p\in M$ and a local normal chart $(U;
x^k)^m_{k=1}$ around it:
$$
v_g(p)=dx^1\wedge\ldots\wedge dx^m, \quad
(i_{X(p)}v_g)_{i\ldots\hat{k}\ldots m}=(-1)^{k+1}\xi^k,
$$
so
$$
\vert\omega\vert^2(p)= \vert
i_{X}v_g\vert^2(p)=\sum_{i_1,\ldots,i_{m-1}=1}^m
\big((i_{X(p)}v_g)_{i_1,\ldots,i_{m-1}}\big)^2=(m-1)! \ \ \vert
X\vert^2(p).
$$
Now, $\int_M\vert X\vert \ v_g= \int_M \vert d\phi.\tau(\phi)\vert
< \infty$, so $\omega$ is integrable. By Lemma
~\ref{eq:lemma-Yau}, $\lim_{i\to\infty}\int_{B_i}d\omega \
v_g=\frac{m-4}{2(m-2)}\lim_{i\to\infty}\int_{B_i}\vert\tau(\phi)\vert^2
\ v_g=0$, hence $\phi$ is harmonic.
\end{proof}

\begin{corollary}
Let $(M,g)$, $m\neq 4$, be an orientable complete Riemannian
manifold and $\phi:(M,g)\to (N,h)$ a map with finite energy and
bienergy. If $S_2=0$ then $\phi$ is harmonic.
\end{corollary}

\begin{proof}
Let $p$ be a point of $M$ and $(U;x^i)^m_{i=1}$,
$(V;y^{\alpha})^n_{\alpha=1}$ be local normal charts around $p$ and
$\phi(p)$, respectively. At $p$
$$
\big\langle d\phi(\frac{\partial}{\partial
x^i}),\tau(\phi)\big\rangle= \big\langle
\phi^{\alpha}_i\frac{\partial}{\partial
y^{\alpha}},\tau(\phi)^{\beta}\frac{\partial}{\partial
y^{\beta}}\big\rangle=
\sum_{\alpha}\phi^{\alpha}_i\tau(\phi)^{\alpha},
$$
so, by Cauchy Inequality
$$
\Big(\big\langle d\phi(\frac{\partial}{\partial
x^i}),\tau(\phi)\big\rangle \Big)^2
\leq\Big(\sum_{\alpha}\big(\phi^{\alpha}_i\big)^2\Big)\vert\tau(\phi)\vert^2.
$$
Therefore, at $p$
$$
\vert d\phi.\tau(\phi)\vert^2=\sum_i\Big(\big\langle
d\phi(\frac{\partial}{\partial x^i}),\tau(\phi)\big\rangle \Big)^2
\leq \vert d\phi\vert^2\vert\tau(\phi)\vert^2.
$$
Consequently,
$$
\int_M\vert d\phi.\tau(\phi)\vert \ v_g \leq \int_M\vert
d\phi\vert\vert\tau(\phi)\vert \ v_g \leq  \Big(\int_M\vert
d\phi\vert^2 \ v_g\Big)^{\frac{1}{2}}
\Big(\int_M\vert\tau(\phi)\vert^2 \ v_g\Big)^{\frac{1}{2}}<\infty,
$$
and we can apply Theorem ~\ref{eq:teorema-Yau}.
\end{proof}

When $m=4$, the situation is drastically different, that is $S_2=0$ does not
always imply $\phi$ harmonic, even if $M$ is
compact.

\begin{theorem}
[\cite{GYJ2}] \label{eq:imersie-pseudo-umbilicala} A non-minimal Riemannian immersion
$\phi:(M^4,g)\to (N,h)$ satisfies
$S_2=0$ if and only if it is pseudo-umbilical.
\end{theorem}

\begin{proof}
First note that for a Riemannian
immersion, $S_2$ reduces to
\begin{equation}
\label{eq:imersie-riemanniana}
S_2(X,Y)=-\frac{1}{2}\vert\tau(\phi)\vert^2\langle X,Y\rangle+
2\langle\tau(\phi),B(X,Y)\rangle,
\end{equation}
$B=\nabla d\phi$ being its second fundamental form.
Recall that a Riemannian immersion is pseudo-umbilical if
and only if its shape operator $A$ satisfies:
$$
A_{\tau(\phi)}=\frac{1}{m}\vert\tau(\phi)\vert^2I,
$$
equivalently
$$
\langle
B(X,Y),\tau(\phi)\rangle=\frac{1}{m}\vert\tau(\phi)\vert^2\langle
X,Y\rangle.
$$
Comparing with~\eqref{eq:imersie-riemanniana} ends the proof.
\end{proof}

To weaken hypotheses, we consider conformal immersions
and show ``rigidity''.

\begin{proposition}
\label{eq:imersie-conforma} Let
$\phi:(M^4,g=e^{2\rho}\phi^{\ast}h)\to (N,h)$ be a conformal
immersion, $M$ compact and orientable. If $S_2=0$ then $\rho$ is
constant and $\ophi:(M^4,\phi^{\ast}h)\to (N,h)$ is a
pseudo-umbilical Riemannian immersion.
\end{proposition}

\begin{proof}
Put $\og=e^{-2\rho}g=\phi^{\ast}h$, denote by $\onabla$
the connection of $\og$, by ${\bf 1}:(M,g)\to
(M,\og)$ the identity map, so that $\ophi:(M,\og)\to (N,h)$ is a
Riemannian immersion, and $\phi=\ophi\circ {\bf 1}$. By chain
rule:
\begin{eqnarray*}
\tau(\phi)&=&d\ophi(\tau({\bf 1}))+\trace_g\nabla d\ophi(d{\bf
1}\cdot,d{\bf 1}\cdot) \\
&=&2d\ophi(\grad_g\rho)+e^{-2\rho}\tau(\ophi),
\end{eqnarray*}
and computing the norm yields
\begin{eqnarray}
\label{eq:conformal-0} \vert\tau(\phi)\vert^2 &=& 4\vert
d\ophi(\grad_g\rho)\vert^2+e^{-4\rho}\vert\tau(\ophi)\vert^2 \\
\nonumber &=& 4e^{-2\rho}\vert\grad_g\rho\vert^2_g
+e^{-4\rho}\vert\tau(\ophi)\vert^2.
\end{eqnarray}
Next:
\begin{eqnarray}
\label{eq:conformal-1} \nonumber
\langle\nabla_Xd\ophi(\grad_g\rho),d\ophi(Y)\rangle&=&
\langle\nabla
d\ophi(X,\grad_g\rho)+d\ophi(\onabla_X\grad_g\rho),d\ophi(Y)\rangle
\\
&=&\og(\onabla_X\grad_g\rho,Y)=\og(\onabla_X(e^{-2\rho}\grad_{\og}\rho),Y)
\\
\nonumber
&=&-2e^{-2\rho}(X\rho)(Y\rho)+e^{-2\rho}\og(\onabla_X\grad_{\og}\rho,Y),
\end{eqnarray}
whilst
\begin{eqnarray}
\label{eq:conformal-2}
\langle\nabla_Xe^{-2\rho}\tau(\ophi),d\phi(Y)\rangle
=-e^{-2\rho}\langle\tau(\ophi),\nabla
d\ophi(X,Y)\rangle.
\end{eqnarray}
From ~\eqref{eq:conformal-1} and ~\eqref{eq:conformal-2}, we
obtain
\begin{eqnarray}
\label{eq:conformal-3} \langle\nabla_X\tau(\phi),d\phi(Y)\rangle
&=&
e^{-2\rho}\big(-4(X\rho)(Y\rho)+2\og(\onabla_X\grad_{\og}\rho,Y)
\\
\nonumber && \ \ \ \ \ \ \ \ -\langle\tau(\ophi),\nabla
d\ophi(X,Y)\rangle\big)
\end{eqnarray}
and
\begin{equation}
\label{eq:conformal-4} \langle d\phi,\nabla\tau(\phi)\rangle=
-e^{-4\rho}\big(4\vert\grad_{\og}\rho\vert^2_{\og}
+2\Delta_{\og}\rho+\vert\tau(\ophi)\vert^2\big),
\end{equation}
where $\Delta_{\og}$ is the Laplacian with respect to $\og$.

\noindent Let $\overline{S_2}$ be the stress-energy
tensor of $\ophi$, by~\eqref{eq:conformal-0},
~\eqref{eq:conformal-3}, ~\eqref{eq:conformal-4}, and
$\vert\grad_g\rho\vert^2_g=e^{-2\rho}\vert\grad_{\og}\rho\vert^2_{\og}$,
a straightforward computation gives:
\begin{eqnarray}
\label{eq:conformal-5} \nonumber
e^{2\rho}S_2(X,Y)&=&\overline{S_2}(X,Y) -2 \big(
\vert\grad_{\og}\rho\vert^2_{\og} + \Delta_{\og}\rho\big)\og(X,Y)
\\
&&+8(X\rho)(Y\rho)-2\big(\og(\onabla_X\grad_{\og}\rho,Y)
+\og(\onabla_Y\grad_{\og}\rho,X)\big).
\end{eqnarray}

Assume $S_2=0$, taking the $\og$-trace of~\eqref{eq:conformal-5},
we obtain $\Delta_{\og}\rho=0$, hence, as $M$ is compact, $\rho$
is constant. Replacing in ~\eqref{eq:conformal-5} shows that
$\overline{S_2}$ vanishes as well.
\end{proof}

\begin{proposition}
Let $\phi:(M^4,g)\to (N,h)$ be a map, $M$ compact, with $S_2=0$.
If there exists $\rho\in C^{\infty}(M)$ such that
$\ophi:(M^4,\og=e^{2\rho}g)\to (N,h)$ is harmonic, then $\phi$ is
harmonic. Moreover, if $\rank \phi=4$, $\rho$ must be constant.
\end{proposition}

\begin{proof}
Denote by ${\bf 1}:(M,g)\to (M,\og)$ the identity map, then
$\tau(\phi)=-2d\ophi(\grad_g\rho)$. Choosing $X=Y=\grad_g\rho$ in
~\eqref{eq:formula-echivalenta}, implies:
\begin{eqnarray*}
0&=&\frac{1}{2}\vert\tau(\phi)\vert^2\vert\grad_g\rho\vert^2_g
+2\langle\nabla_{\grad_g\rho}\tau(\phi),d\phi(\grad_g\rho)\rangle \\
&=& \frac{1}{2}\vert\tau(\phi)\vert^2\vert\grad_g\rho\vert^2_g-
\langle\nabla_{\grad_g\rho}\tau(\phi),\tau(\phi)\rangle \\
&=&\frac{1}{2}\vert\tau(\phi)\vert^2\vert\grad_g\rho\vert^2_g-
\frac{1}{2}\grad_g\rho\big(\vert\tau(\phi)\vert^2\big).
\end{eqnarray*}
Since $M$ is compact, $\vert\tau(\phi)\vert^2$ attains its
maximum at $p_0$. Evaluating the last
equation at this point, shows $\vert\tau(\phi)\vert^2(p_0)=0$,
therefore everywhere.
\newline Furthermore, $0=\tau(\phi)=-2d\phi(\grad_g\rho)$
and $\rank\phi=4$, imply that $\rho$ is constant.
\end{proof}

\begin{proposition}
A map $\phi:(M^4,g)\to (N^4,h)$, $M$ compact, with $S_2=0$ and
$\rank\phi=4$, is harmonic.
\end{proposition}

\begin{proof}
As $d\phi_p:T_pM\to T_{\phi(p)}N$ is an isomorphism at any
point, there exists a unique vector field $Z$ such that
$d\phi_p(Z)=\tau(\phi)(p)$, $\forall p\in M$ and
Equation~\eqref{eq:formula-echivalenta}, with $X=Y=Z$, reads
\begin{eqnarray*}
0&=&\frac{1}{2}\vert\tau(\phi)\vert^2\vert Z\vert^2
+2\langle\nabla_Z\tau(\phi),\tau(\phi)\rangle \\
&=&\frac{1}{2}\vert\tau(\phi)\vert^2\vert Z\vert^2
+Z\big(\vert\tau(\phi)\vert^2\big).
\end{eqnarray*}
Therefore the maximum of $|\tau(\phi)|^2$ must be zero.
\end{proof}

By Proposition~\ref{eq:rang}, if $\phi:(M^4,g)\to (N,h)$ has
$\rank \phi\leq 3$,  $S_2=0$ is equivalent to harmonicity.

\section{Maps with parallel stress-energy tensor}
\label{eq:S-2-paralel}

This section is dedicated to maps with parallel stress-energy
tensor. Before we study in details the condition $\nabla S_2=0$,
we show that, in certain circumstances, this condition is
equivalent to $\Div S_2=0$. Indeed, denoting by $\riem^M$ the
sectional curvature of $M$, we have

\begin{proposition}
Let $\phi:M^m\to N^{m+1}$ be a Riemannian immersion with constant
mean curvature. Assume that $N$ has constant sectional curvature,
and $M$ is compact, orientable with $\riem^M\geq 0$. Then $\Div
S_2=0$ if and only if $\nabla S_2=0$. Moreover, if $\riem^M>0$,
then $\Div S_2=0$ if and only if $S_2=\lambda g$, $\lambda \in\r$,
i.e. $M$ is umbilical.
\end{proposition}

\begin{proof}
Under the hypotheses on $\phi$ and $N$, the Codazzi equation
becomes $(\nabla_XA_H)(Y)=(\nabla_YA_H)(X)$, $\forall X,Y\in
C(TM)$. Thus, the tensor $S_2$, as an $(1,1)$-tensor on $M$,
satisfies
$$
(\nabla_XS_2)(Y)=(\nabla_YS_2)(X).
$$
Now, applying a result of Berger (see, for example, ~\cite[page
202]{PP}), the proposition follows.
\end{proof}

If $\nabla S_2=0$ then clearly $\Div S_2=0$ so, from
Corollary~\ref{eq:corolar-1}, a submersion with parallel
stress-energy tensor is biharmonic. Notice that Proposition
~\ref{eq:aplicatie-S-2} and~\ref{eq:biarmonicitate-identitate}
$a)$, give examples of submersions with $\nabla S_2=0$, and that
Riemannian immersions with $\nabla S_2=0$ have normal bitension
field.

\begin{proposition}
If a non-minimal Riemannian immersion $\phi:(M^m,g)\to (N,h)$,
$m\neq 4$, has $\nabla S_2=0$, then $\vert\tau(\phi)\vert$ is
constant.
\end{proposition}

\begin{proof}
For a Riemannian immersion:
$$
S_2(X,Y)=-\frac{1}{2}\vert\tau(\phi)\vert^2\langle X,Y\rangle+
2\langle\tau(\phi),\nabla d\phi(X,Y)\rangle.
$$
Thus
\begin{eqnarray*}
0&=&(\nabla S_2)(Z,X,Y)=(\nabla_ZS_2)(X,Y) \\
&=&-\frac{1}{2}Z(\vert\tau(\phi)\vert^2)\langle X,Y\rangle
+2\Big( Z\langle\tau(\phi),\nabla d\phi(X,Y)\rangle \\
&&-\langle\tau(\phi),\nabla d\phi(\nabla_ZX,Y)\rangle  -
\langle\tau(\phi),\nabla d\phi(X,\nabla_ZY)\rangle \Big).
\end{eqnarray*}
Take $p\in M$ and $\{X_i\}_{i=1}^m$ a
geodesic frame around it, choose $X=Y=X_i$ and sum up, to obtain:
$$
\frac{4-m}{2}Z(\vert\tau(\phi)\vert^2)=0.
$$
\end{proof}

\begin{remark}
The hypothesis $m\neq 4$ is essential. Indeed, a pseudo-umbilical Riemannian
immersion $\phi:(M^4,g)\to (N,h)$ satisfies $S_2=0$,
hence $\nabla S_2=0$, but $\vert\tau(\phi)\vert$ is not
necessarily constant.
\end{remark}

\begin{theorem}
\label{eq:S-2-paralel-1} A hypersurface ${\bf i}:M^m\to N^{m+1}$,
$m\neq 4$, has $\nabla S_2=0$ if and only if it is parallel.
\end{theorem}

\begin{proof}
If $\nabla S_2=0$, $\vert\tau({\bf i})\vert=m\vert H\vert$ is
constant and $\tau({\bf i})=c\eta$, where $c$ is constant and
$\eta$ is a unit section of the normal bundle. Moreover,
$$
S_2(X,Y)=-\frac{c^2}{2}\langle X,Y\rangle + 2c \ b(X,Y),
$$
where $b(X,Y)=\langle \eta,\nabla d{\bf i}(X,Y)\rangle=\langle
\eta,B(X,Y)\rangle$. We immediately infer that $\nabla S_2=0$
implies $\nabla b=0$, i.e. $M$ is
parallel.

The converse is immediate, since $\nabla b=0$ implies $\vert
H\vert$ constant.
\end{proof}

Parallel hypersurfaces of space forms are classified (\cite{FD}),
and for the Euclidean sphere $\s^{m+1}$ parallel hypersurfaces are
either hyperspheres $\s^m(a)$, $a\in (0,1]$, or Clifford tori
$\s^{m_1}(a_1)\times \s^{m_2}(a_2)$, where $a_1^2+a_2^2=1$,
$a_1\in (0,1)$, and $m_1+m_2=m$. Observe that hyperspheres are
umbilical, while the Clifford tori are not.

\begin{proposition}
\label{eq:S-2-paralel-2} Let $\phi:(M,g)\to (N,h)$ be a non-minimal
pseudo-umbilical Riemannian immersion.
\begin{itemize}
\item[a)] If $m=4$ then $S_2=0$, hence $\nabla S_2=0$.
\item[b)] If $m\neq 4$, $\nabla S_2=0$ if and only if
$\vert\tau(\phi)\vert$ is constant.
\end{itemize}
\end{proposition}

\begin{proof}
Since $\phi$ is pseudo-umbilical, $\langle\nabla
d\phi(X,Y),\tau(\phi)\rangle=\frac{1}{m}\vert\tau(\phi)\vert^2\langle
X,Y\rangle$, consequently,
$$
S_2(X,Y)=\frac{4-m}{2m}\vert\tau(\phi)\vert^2\langle X,Y\rangle.
$$
\end{proof}

The (pseudo-)umbilical hypersurfaces of space forms are classified
\cite{BYC}, and have constant mean curvature. We
indicate now a method to construct pseudo-umbilical submanifolds
of constant mean curvature.

\begin{proposition}
A minimal submanifold ${\bf j}:M^m\to P^{n-1}$ and a non-minimal
umbilical hypersurface of constant mean curvature ${\bf
i}:P^{n-1}\to N^n$ compose into a constant mean curvature
pseudo-umbilical submanifold ${\bf i}\circ {\bf j}:M\to N$.
\end{proposition}

\begin{proof}
Using the chain rule, the proof follows by a standard argument.
\end{proof}

\begin{remark}
Take $N=\s^n$, $P=\s^{n-1}(a)$, $a\in (0,1)$, and $M$ minimal in
$\s^{n-1}(a)$, then the tangent part of $\tau_2({\bf i}\circ {\bf
j})$ vanishes. Moreover, ${\bf i}\circ {\bf j}$ is biharmonic if
and only if $a=\frac{1}{\sqrt{2}}$ (\cite{RCSMCO1,RCSMCO2}).
Differently, replacing $N$ by $\r^n$ or $\h^n$ in the above
construction, ${\bf i}\circ {\bf j}$ is never biharmonic
(\cite{CO}).
\end{remark}

\begin{remark}
Let ${\bf i}:M\to N$ be a submanifold, the tangent
part of $\tau_2({\bf i})$ is (\cite{CO}):
\begin{eqnarray*}
(\tau_2({\bf i}))^{\top}&=& -m\big( \frac{m}{2}\grad (\vert
H\vert^2)+2\trace A_{\nabla^{\perp}_{(\cdot)}}(\cdot)
+2\trace \big( R^N(d{\bf i}\cdot,H)d{\bf i}\cdot\big)^{\top}\big) \\
&=& -m\big(-\frac{m}{2}\grad(\vert H\vert^2) +2\trace (\nabla
A_H)(\cdot,\cdot)\big).
\end{eqnarray*}
A direct computation shows that, indeed, under the hypotheses of
Theorem~\ref{eq:S-2-paralel-1} and Proposition
~\ref{eq:S-2-paralel-2}, $\tau_2({\bf i})$ has
vanishing tangent part.
\end{remark}

\section{The case $S_2=\lambda g$}

The bienergy functional is homogeneous of degree zero with
respect to the domain metric only in dimension four. Therefore, unless
$m=4$, the bienergy can be made arbitrarily large or small by homotheties.
To get around this
problem, one considers variations with fixed volume. 
This type of variational problem is at the heart of Einstein metrics (see
~\cite{MB}).

Let $g$ be a Riemannian metric on a compact, orientable manifold
$M$ and $\{g_t\}$ an isovolumetric variation
of $g$, i.e. $\vol(M,g_t)$ is constant, then
$$
0=\delta(\vol(M,g_t))=\frac{1}{2}\int_M\langle g,\omega\rangle \
v_g.
$$
This says that $\omega$ is orthogonal to $g$ with respect to the
$L^2$-scalar product on $T_gG=C(\odot^2T^{\ast}M)$. By
Theorem~\ref{eq:teorema-prima-variatie}, with respect to
$\{g_t\}$:
$$
\delta(F(g_t))=-\frac{1}{2}\int_M\langle S_2,\omega\rangle \ v_g=
-\frac{1}{2}(S_2,\omega).
$$
Therefore, a critical point $g$ of $F$, with respect to
isovolumetric variations, must be colinear with $S_2$, as vectors in
$T_gG=C(\odot^2T^{\ast}M)$, i.e.
$S_2=\lambda g$ for some $\lambda\in \r$ (\cite{AB,AS}).

The trace of $S_2=\lambda g$ implies:

\begin{proposition}
Let $\phi:(M^m,g)\to (N,h)$ be a map with  $S_2=\lambda g$. If one
of the following holds
\begin{itemize}
\item [a)] $m=2$, \item [b)] $\rank \phi\leq m-1$, $m>2$, \item
[c)] $\phi$ is a submersion, $m>2$,
\end{itemize}
then the norm of $\tau(\phi)$ is constant. Moreover, in case {\rm
c)} $\phi$ is biharmonic.
\end{proposition}

Proposition~\ref{eq:aplicatie-S-2} and
~\ref{eq:biarmonicitate-identitate} $a)$, give examples of
submersions with $S_2=\lambda g$.

\begin{proposition}
If $\phi:(M^m,g)\to (N,h)$ satisfies $S_2=\lambda
g$, with $M$ compact and orientable, then:
\begin{itemize}
\item [a)] If $m=4$, $\lambda=0$.
\item [b)] If $m\neq 4$,
$\frac{4-m}{2}\int_M \vert\tau(\phi)\vert^2 \ v_g=\lambda m\vol (M)$.
\end{itemize}
\end{proposition}

For Riemannian immersions we have the complete classification.

\begin{theorem}
For a non-minimal Riemannian immersion $\phi:(M^m,g)\to (N,h)$,
$S_2=\lambda g$ if and only if $\phi$ is pseudo-umbilical and, if
$m\neq 4$, the norm of $\tau(\phi)$ is constant.
\end{theorem}

\begin{proof}
Since the trace of $S_2=\lambda g$ gives:
\begin{equation}
\label{eq:variatie-izovolumetrica-1}
\frac{4-m}{2}\vert\tau(\phi)\vert^2=\lambda m
\end{equation}
if $m=4$, then $\lambda=0$ and $S_2$ vanishes, i.e.
(Theorem~\ref{eq:imersie-pseudo-umbilicala}) $\phi$ is
pseudo-umbilical.
\newline If $m\neq 4$, ~\eqref{eq:variatie-izovolumetrica-1} shows
that $\vert\tau(\phi)\vert$ is constant. Then, replacing the value
of $\lambda$ in the equation $S_2=\lambda g$, we obtain that $\phi$
is pseudo-umbilical.
\end{proof}

\noindent{\bf Acknowledgements}. The third author wishes to thank
Professor Renzo~Caddeo and the Dipartimento di Matematica e
Informatica, Universit\`a di Cagliari, for hospitality during the
preparation of this paper.

\end{document}